\documentclass[10pt]{amsart}
\usepackage{amsfonts}
\usepackage{color}
\usepackage[colorlinks,
            linkcolor=red,
            anchorcolor=blue,
            citecolor=green
            ]{hyperref}
\usepackage{graphicx}
\usepackage{subfigure}
\usepackage{float}
\usepackage{epstopdf}
\definecolor{c20}{rgb}{0.,0.7,0.}
\definecolor{c30}{rgb}{0.,0.,1.}
\definecolor{c40}{rgb}{1,0.1,0.7}
\definecolor{c50}{rgb}{1,0,0}
\definecolor{c60}{rgb}{1,0.9,0.1}

\def\ehe#1{\textcolor{c30}{#1}}
\def\ehe#1{#1}

\newcommand{\abs}[1]{\left\lvert #1 \right\rvert}

\newcommand{\E}[1]{\mathbb{E}\left\{#1\right\}}

\newcommand{\pk}[1]{\mathbb{P} \left \{#1 \right \} }

\newcommand{\R}{\mathbb{R}}

\newcommand{\BQN}{\begin{eqnarray}}
\newcommand{\EQN}{\end{eqnarray}}
\newcommand{\BQNY}{\begin{eqnarray*}}
\newcommand{\EQNY}{\end{eqnarray*}}

\newcommand{\BS}{\begin{sat}}
\newcommand{\ES}{\end{sat}}
\newcommand{\BT}{\begin{theo}}
\newcommand{\ET}{\end{theo}}
\newcommand{\BL}{\begin{lem}}
\newcommand{\EL}{\end{lem}}
\newcommand{\BK}{\begin{korr}}
\newcommand{\EK}{\end{korr}}

\newcommand{\BD}{\begin{de}}
\newcommand{\ED}{\end{de}}

\newcommand{\BIT}{\begin{itemize}}
\newcommand{\EIT}{\end{itemize}}
\newcommand{\BDI}{\begin{description}}
\newcommand{\EDI}{\end{description}}

\newcommand{\BRM}{\begin{remarks}}
\newcommand{\ERM}{\end{remarks}}

\newcommand{\BEX}{\begin{example}}
\newcommand{\EEX}{\end{example}}

\newtheorem{theo}{Theorem}[section]
\newtheorem{sat}[theo]{Proposition}
\newtheorem{de}[theo]{Definition}
\newtheorem{lem}[theo]{Lemma}

\newtheorem{example}[theo]{Example}
\newtheorem{korr}[theo]{Corollary}

\newtheorem{remarks}[theo]{Remarks}

\newcommand{\prooftheo}[1]{ \textsc{\bf Proof of Theorem} \ref{#1}:}

\newcommand{\prooflem}[1]{\textsc{\bf Proof of Lemma} \ref{#1}:}
\newcommand{\proofkorr}[1]{\textsc{\bf Proof of Corollary} \ref{#1}:}

\newcommand{\COM}[1]{}

\newcommand{\QED}{\hfill $\Box$}

\def\ofp{(\Omega,\mathfrak{F},\mathbb{P})}

\date{}

\def\P{ \mathbb{P}}

\def\E{\mathbb{E}}

\begin{document}

\title[The joint distributions of running maximum of a Slepian processes]{\ehe{The joint distributions of running maximum of a Slepian processes}}

\author{Pingjin Deng}
\address{Pingjin Deng, School of Finance, Nankai University, 300350, Tianjin, PR China,
and Department of Actuarial Science,
University of Lausanne\\
UNIL-Dorigny, 1015 Lausanne, Switzerland
}
\email{Pingjin.Deng@unil.ch}

\maketitle
 {\bf Abstract}:
 Consider the Slepian process $S$ defined by $ S(t)=B(t+1)-B(t),t\in [0,1]$ with
$B(t),t\in \R$ a standard Brownian motion.
  In this contribution we analyze the joint distribution between the  maximum $m_{s}=\max_{0\leq u\leq s}S(u)$  certain and the maximum $M_t=\max_{0\leq u\leq t}S(u)$ for $0< s < t$ fixed.  Explicit integral expression are obtained for the distribution function of the partial maximum $m_{s}$ and the joint distribution function between $m_{s}$ and $M_t$. We also use our results to determine the moments of  $m_{s}$.

\bigskip

\date{\today}
\section{ Introduction} \label{Introduction}
\begin{figure}[H]
 \centering
 \includegraphics[width=0.5\textwidth]{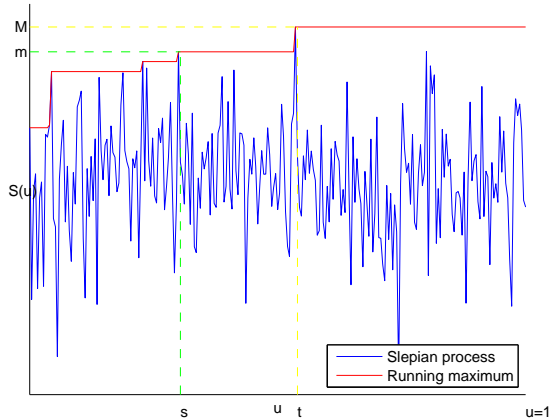}
  \caption{A trajectory of Slepian process (blue line) and its running maximum  (red line) on time interval [0,1]. The partial maxima achieved on time interval $[0,s]$ and  a longer time interval [0,t] are denote by $m$ and $M$. }\label{fig1}
\end{figure}

Throughout this paper, we consider the one-dimensional Slepian process defined as
the increment of a Brownian motion process, namely
\BQN\label{eq1.1}
S(t)=B(t+1)-B(t),\quad t\in [0,1],
\EQN
where $B(t)$ is a standard Brownian motion define on probability space $\ofp$. It can be verified easily that $S(t),\;t\in [0,1]$ is a stationary Gaussian process with covariance function
\BQNY
R_{S}(s,t):=\E[S(s)S(t)]=1-\abs{s-t},\quad s,\;t\in [0,1].
\EQNY
The Slepian processes $S(t)$ which was first defined by Slepian in \cite{slepian1961}, has been studied extensively in stochastic processes and statistics. Zakai and Ziv \cite{zakai1969threshold} gave an application of Slepian processes to the signal shape problem in radar, while the application of these processes to scan statistics and signal dectection problem are presented in  Cressie \cite{cressie1980asymptotic} and Bischoff and Gegg \cite{Bischoff2016}.\\
Another important topic in stochastic processes, where Slepian processes have been wiedly discussed is
the boundary crossing probability. Based on the Markov-like property (or reciprocal property see e.g., \cite{jamison1970reciprocal}) of $S$, Slepian \cite{slepian1961}, Mehr and McFadden \cite{mehr1965certain}, and Shepp \cite{shepp1971first}\cite{shepp1976first} studied the crossing probability of $S$ conditional on $S(0)$ with constant boundary. For a more general boundary, Bischoff and Gegg \cite{Bischoff2016} and Deng \cite{deng2016boundary} gave analytic formulas for the crossing probabilities of $S$ with continuous piecewise linear boundary. For  rencet results on boundary crossing probabilities we refer the reader to \cite{abrahams1984ramp,orsingher1989maximum, 10.2307/2959476, MR2016767,MR2175400, MR2591910, hashorva2014boundary,hashorva2015boundary}.
\newline
For general stochastic processes, both the tail asymptotics of supremum, and the joint survival function of supremum of the process over two intervals has been considered in numerous publications, see e.g.,
 \cite{PickandsB,Berman92, Pit96, HP99, DI2005, DEJ14,Soja,Tabis, Pit20, StamatPopLet, MR3413855, nonhomoANN}. The extremal value statistics are also important in application, for example, the statistics of a maximum is a key process in risk management, the relationship between the risk achieved on a sub-time interval and on the whole time interval can always be characterized using the joint distribution of the running maximum processes. However, the formula of this joint distribution is difficult to establish. In the case of Brownian motion, an explict formula for this joint distribution based on the Fokker-Planck equation is given in \cite{he1998double}.  Recently,the joint distribution between two running maximum both for Brownian motion and Brownian bridge process are studied (see \cite{benichou2016temporal} and \cite{benichou2016joint}, respectively).
\newline
For the Slepian processes defined in equation \eqref{eq1.1}, a little  is known about the partial running maximum and the correlations of different extremes of Slepian process.
This paper is concerned with the maximum statistics of Slepian process $S$.  We obtain an explicit expression for the distribution function of the partial maximum $m_{s}=\max_{0\leq u\leq s}S(u)$. Simple integral expressions are given for this distribution function which allow us to compute the moments generating functions of the running maximum process $m$. We then investigate the joint distribution function between the running maximum $m_{s}$ on a certain time interval $[0,s]$ and $M_t$ on a longer time interval [0,t], see Figure \ref{fig1}. It is interesting that this kind of probability can change into the computation of boundary non-crossing probability of Slepian process with a non-continuous piecewise linear boundary consisting of two lines in finite time interval. Finally, we compute the moments of $m_{s}$ based on its distribution function.
\section{ Results} \label{Results}
In what follows, we let $m_{s}=\max_{0\leq u\leq s}S(u),\;M_t=\max_{0\leq u\leq t}S(u)$, where $S(u)$ is a Slepian process given in \eqref{eq1.1}.
We  aim to compute the following two kinds of probability distribution functions (pdfs): the pdf of the partial maximum $\P(m)$ and  $\P(M)$, the joint distribution  of these two running maxima $\P(m,M)$.
\newline
We start by citing the famous Bachelier-Levy formula  (see e.g. \cite{abundo2016excursions}) which is needed for developing our main results. Concretely, suppose that $a>0$, we have
\BQN\label{eqBL}
\pk{B(t)\leq a+bt,\;\text{for all}\;t\in [0,T]}=\Phi(b\sqrt{T}+\frac{a}{\sqrt{T}})-e^{-2ab}\Phi(b\sqrt{T}-\frac{a}{\sqrt{T}})
\EQN
where $\Phi$ is the distribution of an $N(0,1)$ random variable and the above probability is 0 when $a\leq 0$.
\BRM\label{remark1}
If $b>0,\;T=\infty$, then the probability in equation \eqref{eqBL} is
\BQNY
\pk{B(t)\leq a+bt,\;\text{for all}\;t\geq 0}=1-e^{-2ab}.
\EQNY
\ERM
Next we present our first result for the partial maximum $m_{s}$.
\BT\label{T1}
If $s\in [0,1]$, then the pdf of the running maximum $m_{s}$ of the Slepian process $S$ is given by
\BQN\label{eq2.1}
\P(m)&=&\frac{1}{\sqrt{2\pi}}\int_{-\infty}^{m}\exp\{-\frac{x^2}{2}\}\Phi(\frac{m-x}{2\sqrt{\overline{s}}}+\frac{m+x}{2}\sqrt{\overline{s}})dx\\
\nonumber&&-\frac{1}{\sqrt{2\pi}}\exp\{-\frac{m^2}{2}\}\int_{-\infty}^{m}\Phi(-\frac{m-x}{2\sqrt{\overline{s}}}+\frac{m+x}{2}\sqrt{\overline{s}})dx,
\EQN
where $\overline{s}=\frac{s}{2-s}$.
\ET
The proof of this theorem based on a fact that conditioned on $S(0)$, the Slepian process is equivalent in distribution with a Brownian motion, we give a proof in Section \ref{Proofs}.
\BRM
(i) When $s=0$, then the pdf of $m_0$ is
\BQNY
\P(m)=\frac{1}{\sqrt{2\pi}}\int_{-\infty}^{m}\exp\{-\frac{x^2}{2}\}dx=\Phi(m),
\EQNY
this can also be obtained by the fact that $m_0=S(0)$.\\
(ii) When $s=1$, then from Theorem \ref{T1}, we obtain the  pdf of the global maximum $\max_{0\leq u\leq 1}S(u)$ which we present as follow is also proved in \cite{deng2016boundary},
\BQNY
\P(M)=\pk{\max_{0\leq u\leq 1}S(u)\leq M}=\Phi^2(M)-M\phi(M)\Phi(M)-\phi^2(M),
\EQNY
where $\phi$ is the pdf of $\Phi$; recall $\Phi$ is the df of an $N(0,1)$ random variable.
\ERM
\BRM
If $m=0$ in Theorem \ref{T1}, the probability that the running maximum process $m_s$ take non-positive values is
\BQNY
\P(0)&=&\int_{-\infty}^{0}\phi(x)\Phi(\frac{\overline{s}-1}{2\sqrt{\overline{s}}}x)dx-\frac{1}{\sqrt{2\pi}}\int_{-\infty}^{0}\Phi(\frac{\overline{s}+1}{2\sqrt{\overline{s}}}x)dx\\
&=&\frac{1}{2\pi}\arctan{\frac{2\sqrt{\overline{s}}}{\overline{s}-1}}-\frac{\sqrt{\overline{s}}}{(\overline{s}+1)\pi},
\EQNY
the case $\overline{s}=1$ is Remark 3.2 in \cite{deng2016boundary}.
\ERM
Next, we establish the  joint distribution  function of $m_{s}$ and $M_t$, which is divided into two cases:  $s>0$ and $s=0$.
\BT\label{T2}
If $0 < s\leq t \leq 1$, then the joint pdf of the running maxima $m_{s}$ and $M_t$ of Slepian process $S$ is given by
\BQN\label{eq2.2}
\nonumber\P(m,M)&=&\int_{-\infty}^{m}\int_{-\infty}^{px+q}\frac{1}{2\pi \sqrt{\overline{s}}}\exp\{-\frac{y^2}{2\overline{s}}\}\exp\{-\frac{x^2}{2}\}\Bigl\{1-\exp\{-\frac{(m-x)(px+q-y)}{\overline{s}}\}\Bigr\}\\
&\times& \Bigl\{\Phi(\frac{p x+\eta-y}{\delta}+\frac{M+x}{2}\delta)-\exp\{-(M+x)(p x+\eta-y)\}\Phi(\frac{p x+\eta-y}{\delta}-\frac{M+x}{2}\delta)\Bigr\}dydx,
\EQN
where $p=\frac{1-\overline{s}}{2},\;q=\frac{\overline{s}+1}{2}m,\;\eta=\frac{\overline{s}+1}{2}M,\;\delta=\sqrt{T-\overline{s}},\;\overline{s}=\frac{s}{2-s},\;T=\frac{t}{2-t}$.
\ET
The proof of this theorem is presented in Section \ref{Proofs}.
\BT\label{T3}
If $s=0$, then the joint pdf of the running maxima $m_{0}$ and $M_t$ of Slepian process $S$ is given by
\BQN\label{eq2.3}
\P(m,M)&=&\frac{1}{\sqrt{2\pi}}\int_{-\infty}^{m}\exp\{-\frac{x^2}{2}\}\Phi(\frac{M-x}{2\sqrt{T}}+\frac{M+x}{2}\sqrt{T})dx\\
\nonumber&&-\frac{1}{\sqrt{2\pi}}\exp\{-\frac{M^2}{2}\}\int_{-\infty}^{m}\Phi(-\frac{M-x}{2\sqrt{T}}+\frac{M+x}{2}\sqrt{T})dx,
\EQN
where $T=\frac{t}{2-t}$.
\ET
The proof of this theorem is given in Section \ref{Proofs}.
\subsection{The moments of the partial maximum}
\begin{figure}[H]
\begin{minipage}{0.48\linewidth}
\centerline{\includegraphics[width=0.7\textwidth]{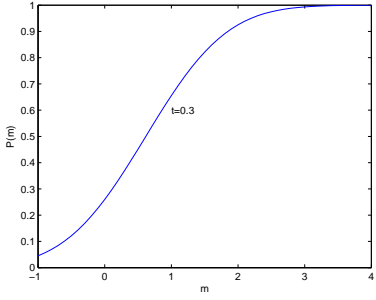}}
\centerline{(a) The distribution of $m_t$ with $t=0.3$.}
\end{minipage}
\hfill
\begin{minipage}{.48\linewidth}
\centerline{\includegraphics[width=0.7\textwidth]{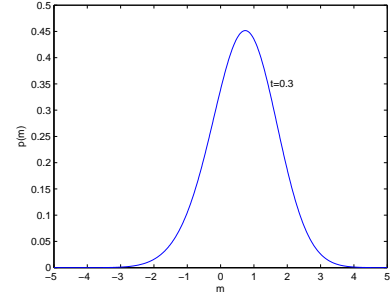}}
\centerline{(b) The density of $m_t$ with $t=0.3$.}
\end{minipage}
\begin{minipage}{0.48\linewidth}
\centerline{\includegraphics[width=0.7\textwidth]{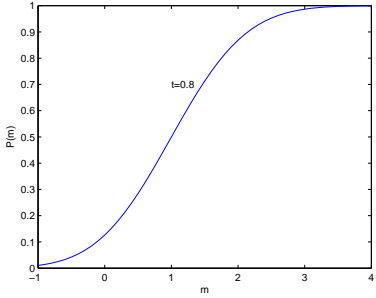}}
\centerline{(c) The distribution of $m_t$ with $t=0.8$.}
\end{minipage}
\hfill
\begin{minipage}{0.48\linewidth}
\centerline{\includegraphics[width=0.7\textwidth]{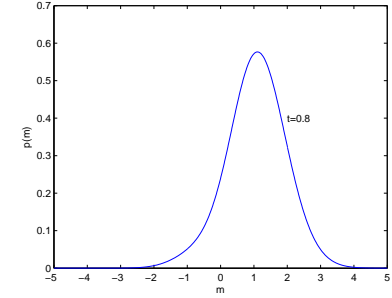}}
\centerline{(d) The density of $m_t$ with $t=0.8$.}
\end{minipage}
\caption{The distribution and density of running maximum $m_t$ given by Equation \eqref{eq2.1} and \eqref{eq2.3} respectively with different $t$.}
\label{fig2}
\end{figure}

Now we begin to compute the moments of the partial maximum $m_{s}$, from Theorem \ref{T1} and after some computation we obtain the density function $p(m)$ of $m_{s}$, which is presented as following:
\BQN\label{eq2.3}
p(m)=\frac{2}{1+\overline{s}}\Phi(\sqrt{\overline{s}}m)\phi(m)
+\frac{2\overline{s}}{1+\overline{s}}m^2\Phi(\sqrt{\overline{s}}m)\phi(m)+\frac{m}{a}\phi(\sqrt{\overline{s}}m)\phi(m),
\EQN
where $a=\frac{1+\overline{s}}{2\sqrt{\overline{s}}}$ is a constant.
From equation \eqref{eq2.3} (or equation \eqref{eq2.1}), we can analysis the features of $m_{s}$. In Figure \ref{fig2}, we plot the distribution and density of running maximum $m_t$.

\begin{figure}[H]
 \centering
\includegraphics[width=0.5\textwidth]{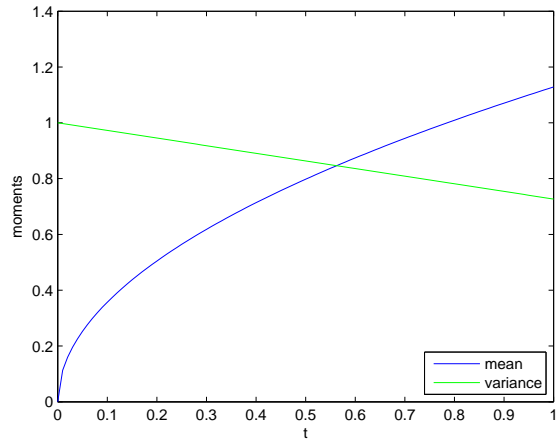}
 \caption{The mean and variance functions of running maximun process $m_t$.}\label{fig3}
\end{figure}

Given $s$, to compute the moments of $m_{s}$,  the moment generating function of $m_{s}$ is given by
\BQNY
M(\theta):=\E[\exp\{\theta m_{s}\}]=\int_{-\infty}^{\infty}\exp\{\theta m\}p(m)dm,
\EQNY
the formula of the k-th moment $\E[m_{s}^k]$  is then given by the k-th derivative of the moment generating function and setting $\theta=0$ , i.e.
\BQNY
\E[m_{s}^k] = \frac{d^{k}M(\theta)}{d\theta^{k}}\left\vert _{\theta=0} \right ..
\EQNY
Using equation \eqref{eq2.3}, we obtain the following:
\BL\label{lemma2.1}
Suppose that $0\leq s\leq 1$ is fixed, the moment generating function of $m_{s}$ is
\BQN\label{eq2.4}
M(\theta)=\exp\{\frac{\theta^2}{2}\}G(\theta),
\EQN
where
\BQNY
G(\theta)=\lambda\int_{-\infty}^{\infty} \Phi(\sqrt{\overline{s}}m)\phi(m-\theta)dm
+\mu\int_{-\infty}^{\infty}m^2\Phi(\sqrt{\overline{s}}m)\phi(m-\theta)dm+\gamma\int_{-\infty}^{\infty}m\phi(\sqrt{\overline{s}}m)\phi(m-\theta)dm,
\EQNY
and
\BQNY
\lambda=\frac{2}{1+\overline{s}},\quad \mu=\frac{2\overline{s}}{1+\overline{s}},\quad \gamma=\frac{2\sqrt{\overline{s}}}{1+\overline{s}},\quad \overline{s}=\frac{s}{2-s}.
\EQNY
\EL
We present the proof of this lemma in Section \ref{Proofs}.  Using equation \eqref{eq2.4}, we can compute the moments for all order, and the first two moments are collected as the following corollary
\BK\label{corollary1}
Given $0\leq s\leq 1$, then the first and second order moments are given by
\BQN\label{eq2.5}
p_1:=\E[m_{s}]=\frac{4\sqrt{\overline{s}}}{\sqrt{2\pi}\sqrt{1+\overline{s}}},
\EQN
\BQN\label{eq2.6}
p_2:=\E[m_{s}^2]=\frac{2+3\overline{s}}{1+\overline{s}}.
\EQN
\EK
The proof of this corollary is displayed in section \ref{Proofs}. Combining equation \eqref{eq2.5} and \eqref{eq2.6}, we can obtain the variance function of $m_{s}$. In Figure \ref{fig3}, we plot The mean and variance functions of $m_t$.

\section{ Proofs} \label{Proofs}
\prooftheo{T1}
Observing that the probability distribution function of the running maximum $m_{s}$ of Slepian processes is
\BQNY
\P(m)&=&\pk{m_{s}=\max_{0\leq u\leq s}S(u)\leq m}\\
&=&\pk{S(u)\leq m,\;\text{for all}\;u\in [0,s]}.
\EQNY
By conditioning on $S(0)$, we represent the above probability as
\BQNY
\P(m)&=&\int_{-\infty}^{m}\pk{S(u)\leq m,\;\text{for all}\;u\in [0,s]\mid S(0)=x}
\varphi(S(0)=x)dx,
\EQNY
where $\varphi(S(0)=x)$ is the density of $S(0)$, i.e.
\BQNY
\varphi(S(0)=x)=\frac{1}{\sqrt{2\pi}}\exp\{-\frac{x^2}{2}\}.
\EQNY
From Lemma 2.3 in \cite{deng2016boundary}, the process $Y=\Bigr\{Y_t=(S(t)\mid S(0)=x),\;t\in [0,1]\Bigl\}$ is equivalent in distribution with process $Z=\Bigr\{Z_t=(2-t)B(\frac{t}{2-t})+(1-t)x,\;t\in [0,1]\Bigl\}$, thus
\BQN\label{eq4.1}
\nonumber &&\P(m)=\int_{-\infty}^{m}\pk{(2-u)B(\frac{u}{2-u})+(1-u)x\leq m,\;\text{for all}\;u\in [0,s]}\varphi(S(0)=x)dx\\
&=&\int_{-\infty}^{m}\pk{B(u)\leq(\frac{m+x}{2})u+\frac{m-x}{2} ,\;\text{for all}\;u\in [0,\frac{s}{2-s}]}\varphi(S(0)=x)dx.
\EQN
Let $\overline{s}=\frac{s}{2-s}$, then from the famous Bachelier-Levy formula (see equation \eqref{eqBL}) we have
\BQN\label{eq4.2}
\nonumber&&\pk{B(u)\leq(\frac{m+x}{2})u+\frac{m-x}{2} ,\;\text{for all}\;u\in [0,\overline{s}]}\\
&=&\Phi(\frac{m-x}{2\sqrt{\overline{s}}}+\frac{m+x}{2}\sqrt{\overline{s}})
-\exp\{-\frac{m^2-x^2}{2}\}\Phi(-\frac{m-x}{2\sqrt{\overline{s}}}+\frac{m+x}{2}\sqrt{\overline{s}}),
\EQN
where $\Phi(x)=\int_{-\infty}^{x}\frac{1}{\sqrt{2\pi}}e^{-\frac{s^2}{2}}ds$ is the cumulative distribution function of standard normal distribution. Substituting equation \eqref{eq4.2} and $\varphi(S(0)=x)=\frac{1}{\sqrt{2\pi}}\exp\{-\frac{x^2}{2}\}$ into equation \eqref{eq4.1}, we conclude that
\BQNY
\P(m)&=&\frac{1}{\sqrt{2\pi}}\int_{-\infty}^{m}\exp\{-\frac{x^2}{2}\}\Phi(\frac{m-x}{2\sqrt{\overline{s}}}+\frac{m+x}{2}\sqrt{\overline{s}})dx\\
\nonumber&&-\frac{1}{\sqrt{2\pi}}\exp\{-\frac{m^2}{2}\}\int_{-\infty}^{m}\Phi(-\frac{m-x}{2\sqrt{\overline{s}}}+\frac{m+x}{2}\sqrt{\overline{s}})dx,
\EQNY
\QED

\prooftheo{T2}
For $0\leq s\leq t\leq 1,\;m\leq M$, the joint probability distribution function between the running maxima $m_{s}$ and $M_t$ of Slepian processes is
\BQNY
\P(m,M)&=&\pk{m_{s}=\max_{0\leq u\leq s}S(u)\leq m,\;M_t=\max_{0\leq u\leq t}S(u)\leq M}\\
&=&\pk{S(u)\leq m,\;\text{for all}\;u\in [0,s]\;\text{and}\;S(u)\leq M \;\text{for all}\;u\in [0,t]}.
\EQNY
Using again the fact that the conditional process $Y=\Bigr\{Y_t=(S(t)\mid S(0)=x),\;t\in [0,1]\Bigl\}$ is equivalent in distribution with process $Z=\Bigr\{Z_t=(2-t)B(\frac{t}{2-t})+(1-t)x,\;t\in [0,1]\Bigl\}$, we obtain
\BQN\label{eq4.3}
\nonumber\P(m,M)&=&\int_{-\infty}^{m}\pk{S(u)\leq m,\;\text{for all}\;u\in [0,s]\;\text{and}\;S(u)\leq M \;\text{for all}\;u\in [0,t]\mid S(0)=x}\varphi(S(0)=x)dx\\
&=&\int_{-\infty}^{m}\P\Bigl\{B(u)\leq(\frac{m+x}{2})u+\frac{m-x}{2} ,\;\text{for all}\;u\in [0,\frac{s}{2-s}]\;\text{and}\;\\
\nonumber &&B(u)\leq(\frac{M+x}{2})u+\frac{M-x}{2} ,\;\text{for all}\;u\in [0,\frac{t}{2-t}]\Bigr\}\varphi(S(0)=x)dx,
\EQN
where $\varphi(S(0)=x)=\frac{1}{\sqrt{2\pi}}\exp\{-\frac{x^2}{2}\}$ is the density of $S(0)$.
Since $0\leq s\leq t\leq 1,\;m\leq M$, then $$\frac{s}{2-s}\leq\frac{t}{2-t};\quad \frac{m+x}{2}u+\frac{m-x}{2}\leq \frac{M+x}{2}u+\frac{M-x}{2},u\in [0,\frac{s}{2-s}],$$ therefore, the last probability in equation \eqref{eq4.3} is equivalent to
\BQNY
\P(m,M)&=&\int_{-\infty}^{m}\P\Bigl\{B(u)\leq(\frac{m+x}{2})u+\frac{m-x}{2} ,\;\text{for all}\;u\in [0,\frac{s}{2-s}]\;\text{and}\;\\
&&B(u)\leq(\frac{M+x}{2})u+\frac{M-x}{2} ,\;\text{for all}\;u\in [\frac{s}{2-s},\frac{t}{2-t}]\Bigr\}\varphi(S(0)=x)dx.
\EQNY
Letting $a=\frac{m+x}{2},\;b=\frac{m-x}{2},\;c=\frac{M+x}{2},\;d=\frac{M-x}{2},\;\overline{s}=\frac{s}{2-s},\;T=\frac{s}{2-s}$, we can simplify $\P(m,M)$ with these notations as
\BQN\label{eq4.4}
\P(m,M)&=&\int_{-\infty}^{m}\P\Bigl\{B(u)\leq au+b ,\;\text{for all}\;u\in [0,\overline{s}]\;\text{and}\;\\
\nonumber&&B(u)\leq cu+d ,\;\text{for all}\;u\in [\overline{s},T]\Bigr\}\varphi(S(0)=x)dx.
\EQN
In fact, denote by
\BQNY
l(u)=
\begin{cases}
au+b,&u\in[0,\overline{s}]\\
cu+d,&u\in[\overline{s},T],
\end{cases}
\EQNY
then equation \eqref{eq4.4} can be viewed as the boundary non-crossing probabilities of Slepian process with piecewise linear function $l(u)$, however, Theorem 3.7 in \cite{deng2016boundary} can not be used here, because $l(u)$ is not continuous at $\overline{s}$. In order to compute $\P(m,M)$ with equation \eqref{eq4.4}, we need compute the non-crossing probabilities of Brownian motion with non-continuous boundary $l(u)$, i.e.
\BQN\label{eq4.5}
\P_{B}^{l}&:=&\P\Bigl\{B(u)\leq au+b ,\;\text{for all}\;u\in [0,\overline{s}]\;\text{and}\;\\
\nonumber&&B(u)\leq cu+d ,\;\text{for all}\;u\in [\overline{s},T]\Bigr\}.
\EQN
The trick here for computing $\P_{B}^{l}$ is using the strong Markovian property of standard Brownian motion $B(u)$ (see e.g.\cite{billingsley2008probability}). Concretely, by conditioning on $B(\overline{s})$ in equation \eqref{eq4.5}, we get
\BQN\label{eq4.6}
\nonumber\P_{B}^{l}&=&\int_{-\infty}^{\min(a\overline{s}+b,c\overline{s}+d)}\P\Bigl\{B(u)\leq au+b ,\;\text{for all}\;u\in [0,\overline{s}]\;\text{and}\;\\
\nonumber&&B(u)\leq cu+d ,\;\text{for all}\;u\in [\overline{s},T]\Big| B(\overline{s})=y\Bigr\}\varphi(B(\overline{s})=y)dy\\
&=&\int_{-\infty}^{a\overline{s}+b}\pk{B(u)\leq au+b ,\;\text{for all}\;u\in [0,\overline{s}]\mid B(\overline{s})=y}\\
\nonumber&&\times \pk{B(u)\leq cu+d ,\;\text{for all}\;u\in [\overline{s},T]\mid B(\overline{s})=y}\varphi(B(\overline{s})=y)dy,
\EQN
where $\varphi(B(\overline{s})=y)=\frac{1}{\sqrt{2\pi \overline{s}}}\exp\{-\frac{y^2}{2\overline{s}}\}$ is the density of $B(\overline{s})$, and the second equality above follows from $a\overline{s}+b\leq c\overline{s}+d$ and the independent property of Brownian motion $B(u)$.
\newline
In equation \eqref{eq4.6},  the rst factor is
\BQNY
&&\pk{B(u)\leq au+b ,\;\text{for all}\;u\in [0,\overline{s}]\mid B(\overline{s})=y}\\
&=&\pk{uB(\frac{1}{u})\leq a+bu ,\;\text{for all}\;u\in [\frac{1}{\overline{s}},\infty)\mid \overline{s}B(\frac{1}{\overline{s}})=y}\\
&=&\pk{B(u)\leq a+bu ,\;\text{for all}\;u\in [\frac{1}{\overline{s}},\infty)\mid B(\frac{1}{\overline{s}})=\frac{y}{\overline{s}}}\\
&=&\pk{B(u)-B(\frac{1}{\overline{s}})\leq a+bu-\frac{y}{\overline{s}} ,\;\text{for all}\;u\in [\frac{1}{\overline{s}},\infty)\mid B(\frac{1}{\overline{s}})=\frac{y}{\overline{s}}}\\
&=&\pk{B(u)-B(\frac{1}{\overline{s}})\leq a+bu-\frac{y}{\overline{s}} ,\;\text{for all}\;u\in [\frac{1}{\overline{s}},\infty)}\\
&=&\pk{B(u)\leq a+b(u+\frac{1}{\overline{s}})-\frac{y}{\overline{s}} ,\;\text{for all}\;u\in [0,\infty)},
\EQNY
the second equality above comes from the fact that $\{uB(\frac{1}{u});\;u\in [\frac{1}{\overline{s}},\infty)\}$ is equivalent in distribution to $\{B(u);\;u\in [0,\overline{s}]\}$, and the last two equalities above hold since the process $\{B(u)-B(\frac{1}{\overline{s}});\;u\in [\frac{1}{\overline{s}},\infty)\}$ is also a standard Brownian motion, and independent with $B(\frac{1}{\overline{s}})$.
From the Bachelier-Levy formula with infinity time horzion (see Remarks \ref{remark1}) we have
\BQNY
\pk{B(u)\leq a+bu+\frac{b}{\overline{s}}-\frac{y}{\overline{s}} ,\;\text{for all}\;u\in [0,\infty)}
=1-\exp\{-\frac{2b(b+a\overline{s}-y)}{\overline{s}}\},
\EQNY
hence the probability
\BQN\label{eq4.7}
\pk{B(u)\leq au+b ,\;\text{for all}\;u\in [0,\overline{s}]\mid B(\overline{s})=y}=1-\exp\{-\frac{2b(b+a\overline{s}-y)}{\overline{s}}\}.
\EQN
Further note that given $B(\overline{s})=y$, the process $B(u+\overline{s})-y$ is again a standard Brownian motion and therefore the second factor in equation \eqref{eq4.6} is
\BQNY
\pk{B(u)\leq cu+d ,\;\text{for all}\;u\in [\overline{s},T]\mid B(\overline{s})=y}
=\pk{B(u)\leq c(u+\overline{s})+d-y,\;\text{for all}\;u\in [0,T-\overline{s}]},
\EQNY
by using the Bachelier-Levy formula again we obtain
\BQN\label{eq4.8}
&&\pk{B(u)\leq cu+d ,\;\text{for all}\;u\in [\overline{s},T]\mid B(\overline{s})=y}\\
\nonumber&=&\Phi(\frac{d+c\overline{s}-y}{\sqrt{T-\overline{s}}}+c\sqrt{T-\overline{s}})-\exp\{-2c(d+c\overline{s}-y)\}\Phi(\frac{d+c\overline{s}-y}{\sqrt{T-\overline{s}}}-c\sqrt{T-\overline{s}}),
\EQN
where $\Phi(x)=\int_{-\infty}^{x}\frac{1}{\sqrt{2\pi}}e^{-\frac{s^2}{2}}ds$ is the cumulative distribution function of standard normal distribution.
\newline
Letting $p=\frac{1-\overline{s}}{2},\;q=\frac{\overline{s}+1}{2}m,\;\eta=\frac{\overline{s}+1}{2}M,\;
\delta=\sqrt{T-\overline{s}},\;\overline{s}=\frac{s}{2-s},\;T=\frac{t}{2-t}$, and substituting equation \eqref{eq4.7} and equation \eqref{eq4.8} into equation \eqref{eq4.6} we conclude that
\BQNY
\P(m,M)&=&\int_{-\infty}^{m}\int_{-\infty}^{px+q}\frac{1}{2\pi \sqrt{\overline{s}}}\exp\{-\frac{y^2}{2\overline{s}}\}\exp\{-\frac{x^2}{2}\}\Bigl\{1-\exp\{-\frac{(m-x)(px+q-y)}{\overline{s}}\}\Bigr\}\\
&\times& \Bigl\{\Phi(\frac{p x+\eta-y}{\delta}+\frac{M+x}{2}\delta)-\exp\{-(M+x)(p x+\eta-y)\}\Phi(\frac{p x+\eta-y}{\delta}-\frac{M+x}{2}\delta)\Bigr\}dydx,
\EQNY
completing the proof.
\QED

\prooftheo{T3}
For $s=0,\;m\leq M$, the joint probability distribution function between the running maxima $m_{0}$ and $M_t$ of Slepian processes is
\BQNY
\P(m,M)&=&\pk{m_{0}=S(0)\leq m,\;M_t=\max_{0\leq u\leq t}S(u)\leq M}\\
&=&\pk{S(0)\leq m,\;\text{and}\;S(u)\leq M \;\text{for all}\;u\in [0,t]}.
\EQNY
Conditioning on $S(0)$ and using the same method as in the proof of Theorem \ref{T1}, we have
\BQNY
\P(m,M)&=&\frac{1}{\sqrt{2\pi}}\int_{-\infty}^{m}\exp\{-\frac{x^2}{2}\}\Phi(\frac{M-x}{2\sqrt{T}}+\frac{M+x}{2}\sqrt{T})dx\\
&&-\frac{1}{\sqrt{2\pi}}\exp\{-\frac{M^2}{2}\}\int_{-\infty}^{m}\Phi(-\frac{M-x}{2\sqrt{T}}+\frac{M+x}{2}\sqrt{T})dx,
\EQNY
where $T=\frac{t}{2-t}$, then the claim follows. 
\QED

\prooflem{lemma2.1}
Since
\BQNY
M(\theta):=\E[\exp\{\theta m_{s}\}],
\EQNY
from equation \eqref{eq2.3} we obtain
\BQNY
M(\theta)=\int_{-\infty}^{\infty}\exp\{\theta m\}\Bigl\{\frac{2}{1+\overline{s}}\Phi(\sqrt{\overline{s}}m)\phi(m)
+\frac{2\overline{s}}{1+\overline{s}}m^2\Phi(\sqrt{\overline{s}}m)\phi(m)+\frac{m}{a}\phi(\sqrt{\overline{s}}m)\phi(m)\Bigr\}dm,
\EQNY
where $a=\frac{1+\overline{s}}{2\sqrt{\overline{s}}}$. Let $ \lambda=\frac{2}{1+\overline{s}},\; \mu=\frac{2\overline{s}}{1+\overline{s}},\; \gamma=\frac{2\sqrt{\overline{s}}}{1+\overline{s}},$ then we have
\BQNY
M(\theta)=\lambda\int_{-\infty}^{\infty}\exp\{\theta m\}\Phi(\sqrt{\overline{s}}m)\phi(m)dm+
\mu\int_{-\infty}^{\infty}m^2\exp\{\theta m\}\Phi(\sqrt{\overline{s}}m)\phi(m)dm+\gamma\int_{-\infty}^{\infty}m\exp\{\theta m\}\phi(\sqrt{\overline{s}}m)\phi(m)\Bigr\}dm.
\EQNY
Observing that for any $\theta\in \R$, we have
\BQN\label{eq4.9}
\nonumber \exp\{\theta m\}\Phi(\sqrt{\overline{s}}m)\phi(m)&=&\frac{1}{\sqrt{2\pi}}\exp\{\theta m\}\exp\{-\frac{m^2}{2}\}\Phi(\sqrt{\overline{s}}m)\\
\nonumber &=&\frac{1}{\sqrt{2\pi}}\exp\{-\frac{1}{2}(m-\theta)^2\}\exp\{\frac{\theta^2}{2}\}\Phi(\sqrt{\overline{s}}m)\\
&=&\exp\{\frac{\theta^2}{2}\}\Phi(\sqrt{\overline{s}}m)\phi(m-\theta)
\EQN
similarly, we have
\BQN\label{eq4.10}
\nonumber \exp\{\theta m\}\phi(\sqrt{\overline{s}}m)\phi(m)&=&\frac{1}{\sqrt{2\pi}}\exp\{\theta m\}\exp\{-\frac{m^2}{2}\}\phi(\sqrt{\overline{s}}m)\\
&=&\exp\{\frac{\theta^2}{2}\}\phi(\sqrt{\overline{s}}m)\phi(m-\theta).
\EQN
Substituting equation \eqref{eq4.9} and \eqref{eq4.10} into $M(\theta)$, and let
\BQNY
G(\theta)=\lambda\int_{-\infty}^{\infty} \Phi(\sqrt{\overline{s}}m)\phi(m-\theta)dm
+\mu\int_{-\infty}^{\infty}m^2\Phi(\sqrt{\overline{s}}m)\phi(m-\theta)dm+\gamma\int_{-\infty}^{\infty}m\phi(\sqrt{\overline{s}}m)\phi(m-\theta)dm,
\EQNY
then the lemma established.
\QED

\proofkorr{corollary1}
Taking the first derivative of equation \eqref{eq2.4} and letting $\theta=0$, we have
\BQNY
p_1=\lambda\int_{-\infty}^{\infty} m\Phi(\sqrt{\overline{s}}m)\phi(m)dm
+\mu\int_{-\infty}^{\infty}m^3\Phi(\sqrt{\overline{s}}m)\phi(m)dm+\gamma\int_{-\infty}^{\infty}m^2\phi(\sqrt{\overline{s}}m)\phi(m)dm,
\EQNY
where $ \lambda=\frac{2}{1+\overline{s}},\; \mu=\frac{2\overline{s}}{1+\overline{s}},\; \gamma=\frac{2\sqrt{\overline{s}}}{1+\overline{s}}$. It is easily to check that
\BQNY
\int m\Phi(\sqrt{\overline{s}}m)\phi(m)dm&=&\frac{\sqrt{\overline{s}}}{\sqrt{2\pi}\sqrt{1+\overline{s}}}\Phi(m\sqrt{1+\overline{s}})-\Phi(\sqrt{\overline{s}}m)\phi(m)+C_1,\\
\int m^3\Phi(\sqrt{\overline{s}}m)\phi(m)dm&=&\frac{2T_{1}^{\frac{3}{2}}+3\sqrt{\overline{s}}}{\sqrt{2\pi}(1+\overline{s})^{\frac{3}{2}}}\Phi(m\sqrt{1+\overline{s}})-(m^2+2)\Phi(\sqrt{\overline{s}}m)\phi(m)-\frac{\sqrt{\overline{s}}m}{\sqrt{2\pi}(1+\overline{s})}\phi(\sqrt{\overline{s}}m)+C_2,
\EQNY
where $C_1,\;C_2$ are constant. Thus, we have
\BQNY
a_1&:=&\int_{-\infty}^{\infty} m\Phi(\sqrt{\overline{s}}m)\phi(m)dm=\frac{\sqrt{\overline{s}}}{\sqrt{2\pi}\sqrt{1+\overline{s}}},\\
a_2&:=&\int_{-\infty}^{\infty}m^3\Phi(\sqrt{\overline{s}}m)\phi(m)dm=\frac{2T_{1}^{\frac{3}{2}}+3\sqrt{\overline{s}}}{\sqrt{2\pi}(1+\overline{s})^{\frac{3}{2}}}.
\EQNY
Using the integral by part formula, we have
\BQNY
a_3&:=&\int_{-\infty}^{\infty}m^2\phi(\sqrt{\overline{s}}m)\phi(m)dm\\
&=&\frac{1}{\sqrt{\overline{s}}}\{\int_{-\infty}^{\infty}m^3\Phi(\sqrt{\overline{s}}m)\phi(m)dm-2\int_{-\infty}^{\infty} m\Phi(\sqrt{\overline{s}}m)\phi(m)dm\}\\
&=&\frac{1}{\sqrt{\overline{s}}}(a_2-2a_1),
\EQNY
Hence we obtain that
\BQNY
p_1=\lambda a_1+\mu a_2+\gamma a_3=\frac{4\sqrt{\overline{s}}}{\sqrt{2\pi}\sqrt{1+\overline{s}}}.
\EQNY
Taking second derivative of equation \eqref{eq2.4} and letting $\theta=0$, we have
\BQNY
p_2=\lambda\int_{-\infty}^{\infty} m^2\Phi(\sqrt{\overline{s}}m)\phi(m)dm
+\mu\int_{-\infty}^{\infty}m^4\Phi(\sqrt{\overline{s}}m)\phi(m)dm+\gamma\int_{-\infty}^{\infty}m^3\phi(\sqrt{\overline{s}}m)\phi(m)dm,
\EQNY
by an analogy method we get
\BQNY
&&\int_{-\infty}^{\infty} m^2\Phi(\sqrt{\overline{s}}m)\phi(m)dm=\frac{1}{2},\\
&&\int_{-\infty}^{\infty}m^4\Phi(\sqrt{\overline{s}}m)\phi(m)dm=\frac{3}{2},\\
&&\int_{-\infty}^{\infty}m^2\phi(\sqrt{\overline{s}}m)\phi(m)dm=0.
\EQNY
Hence we have
\BQNY
p_2:=\E[m_{s}^2]=\frac{2+3\overline{s}}{1+\overline{s}}
\EQNY
establishing the proof.
\QED
\section{Acknowledgement}
This work was partly financed by the project NSFC Grant NO.71573143 and SNSF Grant 200021-166274.

\bibliographystyle{ieeetr}
\bibliography{RunningE}

\def\polhk#1{\setbox0=\hbox{#1}{\ooalign{\hidewidth
  \lower1.5ex\hbox{`}\hidewidth\crcr\unhbox0}}}
\begin{thebibliography}{10}

\bibitem{slepian1961}
D.~Slepian, ``First passage time for a particular gaussian process,'' {\em Ann.
  Math. Statist.}, vol.~32, pp.~610--612, 06 1961.

\bibitem{zakai1969threshold}
M.~Zakai and J.~Ziv, ``On the threshold effect in radar range estimation
  (corresp.),'' {\em IEEE Transactions on Information Theory}, vol.~15, no.~1,
  pp.~167--170, 1969.

\bibitem{cressie1980asymptotic}
N.~Cressie, ``The asymptotic distribution of the scan statistic under
  uniformity,'' {\em The Annals of Probability}, pp.~828--840, 1980.

\bibitem{Bischoff2016}
W.~Bischoff and A.~Gegg, ``Boundary crossing probabilities for
  (q,d)-slepian-processes,'' {\em Statistics and Probability Letters},
  pp.~1--6, 2016.

\bibitem{jamison1970reciprocal}
B.~Jamison, ``Reciprocal processes: The stationary gaussian case,'' {\em The
  Annals of Mathematical Statistics}, vol.~41, no.~5, pp.~1624--1630, 1970.

\bibitem{mehr1965certain}
C.~Mehr and J.~McFadden, ``Certain properties of gaussian processes and their
  first-passage times,'' {\em Journal of the Royal Statistical Society. Series
  B (Methodological)}, pp.~505--522, 1965.

\bibitem{shepp1971first}
L.~Shepp, ``First passage time for a particular gaussian process,'' {\em The
  Annals of Mathematical Statistics}, pp.~946--951, 1971.

\bibitem{shepp1976first}
L.~Shepp and D.~Slepian, ``First-passage time for a particular stationary
  periodic gaussian process,'' {\em Journal of Applied Probability},
  pp.~27--38, 1976.

\bibitem{deng2016boundary}
P.~Deng, ``The boundary non-crossing probabilities for slepian process,'' {\em
  arXiv preprint arXiv:1608.01133}, 2016.

\bibitem{abrahams1984ramp}
J.~Abrahams, ``Ramp crossings for slepian's process,'' {\em IEEE transactions
  on information theory}, vol.~30, no.~3, pp.~574--575, 1984.

\bibitem{orsingher1989maximum}
E.~Orsingher, ``On the maximum of gaussian fourier series emerging in the
  analysis of random vibrations,'' {\em Journal of Applied Probability},
  pp.~182--188, 1989.

\bibitem{10.2307/2959476}
I.~B.-D. Moshe Ein-Gal, ``Passages and maxima for a particular gaussian
  process,'' {\em The Annals of Probability}, vol.~3, no.~3, pp.~549--556,
  1975.

\bibitem{MR2016767}
W.~Bischoff, F.~Miller, E.~Hashorva, and J.~H{\"u}sler, ``Asymptotics of a
  boundary crossing probability of a {B}rownian bridge with general trend,''
  {\em Methodol. Comput. Appl. Probab.}, vol.~5, no.~3, pp.~271--287, 2003.

\bibitem{MR2175400}
E.~Hashorva, ``Exact asymptotics for boundary crossing probabilities of
  {B}rownian motion with piecewise linear trend,'' {\em Electron. Comm.
  Probab.}, vol.~10, pp.~207--217 (electronic), 2005.

\bibitem{MR2591910}
E.~Hashorva, ``Boundary non-crossings of {B}rownian pillow,'' {\em J. Theoret.
  Probab.}, vol.~23, no.~1, pp.~193--208, 2010.

\bibitem{hashorva2014boundary}
E.~Hashorva and Y.~Mishura, ``Boundary noncrossings of additive wiener
  fields?,'' {\em Lithuanian Mathematical Journal}, vol.~54, no.~3,
  pp.~277--289, 2014.

\bibitem{hashorva2015boundary}
E.~Hashorva, Y.~Mishura, and O.~Seleznjev, ``Boundary non-crossing
  probabilities for fractional {B}rownian motion with trend,'' {\em Stochastics
  An International Journal of Probability and Stochastic Processes}, vol.~87,
  no.~6, pp.~946--965, 2015.

\bibitem{PickandsB}
J.~Pickands, III, ``Asymptotic properties of the maximum in a stationary
  {G}aussian process,'' {\em Trans. Amer. Math. Soc.}, vol.~145, pp.~75--86,
  1969.

\bibitem{Berman92}
S.~Berman, {\em Sojourns and extremes of stochastic processes}.
\newblock The Wadsworth \& Brooks/Cole Statistics/Probability Series, Pacific
  Grove, CA: Wadsworth \& Brooks/Cole Advanced Books \& Software, 1992.

\bibitem{Pit96}
V.~I. Piterbarg, {\em Asymptotic methods in the theory of {G}aussian processes
  and fields}, vol.~148 of {\em Translations of Mathematical Monographs}.
\newblock Providence, RI: American Mathematical Society, 1996.

\bibitem{HP99}
J.~H{\"u}sler and V.~Piterbarg, ``Extremes of a certain class of {G}aussian
  processes,'' {\em Stochastic Process. Appl.}, vol.~83, no.~2, pp.~257--271,
  1999.

\bibitem{DI2005}
A.~Dieker, ``Extremes of {G}aussian processes over an infinite horizon,'' {\em
  Stochastic Process. Appl.}, vol.~115, no.~2, pp.~207--248, 2005.

\bibitem{DEJ14}
K.~D{\c{e}}bicki, E.~Hashorva, and L.~Ji, ``Tail asymptotics of supremum of
  certain {G}aussian processes over threshold dependent random intervals,''
  {\em Extremes}, vol.~17, no.~3, pp.~411--429, 2014.

\bibitem{Soja}
K.~D{\c{e}}bicki, E.~Hashorva, and N.~Soja-Kukie{\l}a, ``Extremes of
  homogeneous {G}aussian random fields,'' {\em J. Appl. Probab.}, vol.~52,
  no.~1, pp.~55--67, 2015.

\bibitem{Tabis}
K.~D{\c{e}}bicki, E.~Hashorva, L.~Ji, and K.~Tabi{\'s}, ``Extremes of
  vector-valued {G}aussian processes: {E}xact asymptotics,'' {\em Stochastic
  Process. Appl.}, vol.~125, no.~11, pp.~4039--4065, 2015.

\bibitem{Pit20}
V.~I. Piterbarg, {\em Twenty Lectures About {G}aussian Processes}.
\newblock London, New York: Atlantic Financial Press, 2015.

\bibitem{StamatPopLet}
G.~Popivoda and S.~Stamatovic, ``Extremes of {G}aussian fields with a smooth
  random variance,'' {\em Statist. Probab. Lett.}, vol.~110, pp.~185--190,
  2016.

\bibitem{MR3413855}
E.~Hashorva and L.~Ji, ``Extremes of {$\alpha(\bold{t})$}-locally stationary
  {G}aussian random fields,'' {\em Trans. Amer. Math. Soc.}, vol.~368, no.~1,
  pp.~1--26, 2016.

\bibitem{nonhomoANN}
K.~D\c{e}bicki, E.~Hashorva, and L.~Ji, ``Extremes of a class of nonhomogeneous
  gaussian random fields,'' {\em Ann. Probab.}, vol.~44, no.~2, pp.~984--1012,
  2016.

\bibitem{he1998double}
H.~He, W.~P. Keirstead, and J.~Rebholz, ``Double lookbacks,'' {\em Mathematical
  Finance}, vol.~8, no.~3, pp.~201--228, 1998.

\bibitem{benichou2016temporal}
O.~Benichou, P.~Krapivsky, C.~Mejia-Monasterio, and G.~Oshanin, ``Temporal
  correlations of the running maximum of a brownian trajectory,'' {\em arXiv
  preprint arXiv:1602.06770}, 2016.

\bibitem{benichou2016joint}
O.~B{\'e}nichou, P.~Krapivsky, C.~Mej{\'\i}a-Monasterio, and G.~Oshanin,
  ``Joint distributions of partial and global maxima of a brownian bridge,''
  {\em Journal of Physics A: Mathematical and Theoretical}, vol.~49, no.~33,
  p.~335002, 2016.

\bibitem{abundo2016excursions}
M.~Abundo, ``On the excursions of drifted brownian motion and the successive
  passage times of brownian motion,'' {\em Physica A: Statistical Mechanics and
  its Applications}, vol.~457, pp.~176--182, 2016.

\bibitem{billingsley2008probability}
P.~Billingsley, {\em Probability and measure}.
\newblock John Wiley \& Sons, 2008.

\end{thebibliography}
\end{document}